\documentclass{amsart}

\usepackage[usenames]{color}
\usepackage{amssymb}
\usepackage{graphicx, epsfig}
\usepackage{latexsym, amsfonts, amscd, amsmath}
\usepackage{mathrsfs}
\makeindex 
\input xy
\xyoption{all}

\definecolor{Indigo}{rgb}{0.2,0.1,0.7}
\definecolor{Violet}{rgb}{0.5,0.1,0.7}

\newtheorem{thm}{Theorem}[section]
\newtheorem{prop}[thm]{Proposition}
\newtheorem{lem}[thm]{Lemma}
\newtheorem{cor}[thm]{Corollary}

\theoremstyle{definition}
\newtheorem{dfn}[thm]{Definition}

\theoremstyle{remark}
\newtheorem{rmk}[thm]{Remark}

\newenvironment{mat}{\left(\begin{smallmatrix}
\end{smallmatrix}\right)}{enddef}

\numberwithin{equation}{section} \numberwithin{figure}{section}
\numberwithin{table}{section}

\newcommand{\Aut}{{\operatorname{Aut}}}

\newcommand{\End}{{\operatorname{End}}}
\newcommand{\Fr}{{\operatorname{Fr }}}

\newcommand{\Ker}{{\operatorname{Ker}}}

\newcommand{\Spec}{{\operatorname{Spec }}}
\newcommand{\Spf}{{\operatorname{Spf }}}
\newcommand{\Sp}{{\operatorname{Spm }}}

\newcommand{\val}{{\operatorname{val}}}
\newcommand{\Ver}{{\operatorname{Ver }}}

\newcommand{\GL}{{\operatorname{GL}}}

\newcommand{\gera}{{\frak{a}}}

\newcommand{\germ}{{\frak{m}}}

\newcommand{\gerp}{{\frak{p}}}

\newcommand{\gers}{{\frak{s}}}
\newcommand{\gert}{{\frak{t}}}

\newcommand{\gerB}{{\frak{B}}}

\newcommand{\gerI}{{\frak{I}}}

\newcommand{\gerX}{{\frak{X}}}
\newcommand{\gerY}{{\frak{Y}}}
\newcommand{\gerZ}{{\frak{Z}}}

\newcommand{\calO}{{\mathcal{O}}}

\newcommand{\calU}{{\mathcal{U}}}
\newcommand{\calV}{{\mathcal{V}}}

\newcommand{\calX}{{\mathcal{X}}}

\def\AA{\mathbb{A}}

\def\CC{\mathbb{C}}

\def\FF{\mathbb{F}}

\def\HH{\mathbb{H}}

\def\NN{\mathbb{N}}

\def\PP{\mathbb{P}}
\def\QQ{\mathbb{Q}}
\def\RR{\mathbb{R}}

\def\ZZ{\mathbb{Z}}

\newcommand{\scrA}{{\mathscr{A}}}

\newcommand{\scrE}{{\mathscr{E}}}

\newcommand{\scrG}{{\mathscr{G}}}
\newcommand{\scrH}{{\mathscr{H}}}
\newcommand{\scrI}{{\mathscr{I}}}

\newcommand{\scrP}{{\mathscr{P}}}

\newcommand{\scrZ}{{\mathscr{Z}}}

\newcommand{\id}{{\noindent}}

\newcommand{\arr}{{\; \longrightarrow \;}}

\newcommand{\injects}{{\; \hookrightarrow \;}}

\newcommand{\ol}{{\mathcal{O}_L}}

\newcommand{\Xrig}{{\gerX_\text{\rm rig}}}

\newcommand{\Yrig}{{\gerY_\text{\rm rig}}}

\newcommand{\pirig}{{\pi_\text{\rm rig}}}
\newcommand{\rig}{{\operatorname{rig}}}
\newcommand{\spe}{{\operatorname{sp }}}

\begin{document}
\marginparwidth 50pt

\title[Canonical Subgroup]{The canonical subgroup: a ``subgroup-free" approach}
\author{Eyal Z. Goren \& Payman L Kassaei}
\email{goren@math.mcgill.ca; kassaei@math.mcgill.ca}
\address{Department of Mathematics and Statistics,
McGill University, 805 Sherbrooke West, Montreal H3A 2K6, Quebec,
Canada}

\subjclass{Primary [11F85, 11F33]; Secondary [11G18, 14G35,
14G22]}

\keywords{Canonical subgroup, overconvergent modular form, rigid
geometry.}

\date{January 7, 2005
}
\maketitle

\section{Introduction}
\id Canonical subgroups are essential to the theory of
overconvergent modular forms. An elliptic curve~$E$ with an
ordinary reduction modulo a prime~$p$ has a distinguished subgroup
of rank~$p$, which is the kernel of multiplication by~$p$ on its
formal group. This subgroup is a canonical lift of the kernel of
$\Fr_p$ on~$E$ modulo~$p$. The overconvergence of the canonical
subgroup, i.e. the fact that it can also be defined for elliptic
curves with a ``not too supersingular'' reduction modulo~$p$,
allows one to define and study the~${\mathbf U}_p$ operator for
overconvergent modular forms (See \cite[\S3.11]{Katz}.). Recently,
in \cite{Buzzard,Kassaei3}, this theory has found new applications
to the problem of analytic continuation of overconvergent modular
forms. In these applications it is essential to understand the
\emph{precise extent} of overconvergence of the canonical
subgroup, and to determine the ``measure of supersingularity" of a
quotient of an elliptic curve by a subgroup of order~$p$
(including the canonical subgroup). These results appear in
\cite[Thms. 3.1, 3.10.7]{Katz}, where they are attributed to
Lubin. A slightly more general version can be found in
\cite{Buzzard}.

Classically, the canonical subgroup of an elliptic curve (when it
exists) is constructed by a close study of the power series of
multiplication by~$p$ in its formal group. In
\cite{Kassaei1,Kassaei2} this approach was used to develop a
similar theory over certain PEL Shimura curves. Generalizing this
approach to higher dimensions seems to pose a serious challenge,
because it uses the one-dimensionality of the formal group in an
essential way, including the existence of Newton polygon for power
series in one variable.

The problem of constructing a canonical subgroup for each elliptic
curve belonging to a certain region of a modular curve~$X(\Gamma)$
can be rephrased as finding a partial section to the forgetful
morphism of rigid analytic curves~$\pi\colon X(\Gamma_0(p) \cap
\Gamma) \rightarrow X(\Gamma)$ whose moduli-theoretic description
is $(E,\gamma,H) \mapsto (E,\gamma)$ where $(E,\gamma)$ is an
elliptic curve with level $\Gamma$-structure and $H \subset E[p]$
is a finite flat subgroup of order $p$. Our approach ignores this
moduli-theoretic description and just takes into account the
geometry of the morphism $\pi$. This is what we call the
``subgroup-free" approach. It has been known for a while that one
can prove, using a general principle of rigid geometry due to
Berthelot \cite{Berthelot}, that such a section defined over the
ordinary locus overconverges (to an a priori non-tractable extent)
beyond the ordinary locus. This approach, which is expected to
work in other situations, was used in \cite{KisinLai} to prove the
overconvergence of canonical subgroups in the case of Hilbert
modular varieties. However, other aspects of the theory, which
were discussed in the opening paragraph, remain unsettled even in
the case of Hilbert modular varieties. These aspects are also not
fully covered by other recent approaches \cite{AbbesMokrane,
AndreattaGasbarri, Nevens}.

The purpose of this article is to derive all aspects of the theory
of canonical subgroups via  the``subgroup-free"  approach. Our
thesis is that the rigid geometric (or formal schematic) picture
that arises in the familiar setting of the relevant Shimura
varieties suffices by itself to guarantee the existence of the
canonical subgroup and many of its properties. In this manuscript
we demonstrate that for Shimura varieties of dimension one, even
if they do not possess a natural modular interpretation. In fact,
this lack of a moduli interpretation can be taken as a further
motivation for our approach. Notice that our approach is such that
inspires generalization to higher dimensional settings. More
specifically, one constructs a section over the ordinary locus by
lifting a section from characteristic $p$. One separately studies
sections over the non-ordinary locus by using the theory of local
models for the special fibre of the Shimura variety in question,
and finally these two sections are glued together by using the
above-mentioned principle of rigid geometry along with a certain
uniqueness result. The authors hope to pursue this subject in a
future publication.

\

\id Let~$p$ be a prime. Let~$\calO=\calO_K$ be the ring of
integers of a finite extension~$K$ of~$\QQ_p$, $\varpi$ a
uniformizer of~$\calO$, $\kappa = \calO/(\varpi)$ the residue
field, and~$\val=\val_K$ be the valuation normalized so that
$\val(\varpi)=1$. By a ``curve"~$X$ over~$\calO$ we mean a flat finite-type
morphism~$f\colon X \rightarrow \calO$ of relative dimension~$1$
of a reduced separated scheme~$X$, such that the geometric fibres
of~$f$ are connected;~$f$ need not be proper.

Let~$X, Y$ be curves over~$\calO$. We assume that~$X, Y$ are
regular schemes, $X \rightarrow \Spec(\calO)$ is smooth
and~$\pi\colon Y \arr X$ is a finite flat morphism of
degree~$e+1$. Moreover, we assume that (i) there exists a
section~$s\colon X\otimes \kappa \arr Y\otimes \kappa$ to~$\pi
\otimes \kappa$, that (ii) the special fibre~$Y\otimes \kappa$ is
a reduced normal crossing divisor with two components, and (iii)
the set theoretic preimage~$(\pi\otimes \kappa)^{-1}(\pi\otimes
\kappa)(Q)$ is equal to~$Q$ for any singular point~$Q\in Y\otimes
\kappa$. To remove any doubt, we assume that~$Y \otimes \kappa$
\underline{is} singular and by a normal crossing divisor we mean
that each intersection is defined over~$\kappa$ and its completed
local ring is isomorphic to~$\kappa[\![s, t]\!]/(st)$. We
define~$(Y \otimes \kappa)^\infty=s(X \otimes \kappa)\setminus (Y
\otimes \kappa)^{\rm sing}$, and~$(Y \otimes \kappa)^0=(Y \otimes
\kappa) \setminus s(X \otimes \kappa)$.

From the point of view of a general theory this is a very specific
situation, nonetheless it (and its appropriate generalization) is
the one that occurs for Shimura curves (respectively,
higher-dimensional PEL Shimura varieties); see
\S\ref{section:applications}. In fact, condition (iii) is only put
to have ``cleaner statements"; it holds in the case of Shimura
curves. Under these conditions, we prove in \S\ref{section:main
theorem} the following result.

Let~$\gerX, \gerY$ be the formal schemes obtained, respectively,
by completing~$X, Y$ along their special fibres. The induced
morphism~$\gerY \rightarrow \gerX$ is still denoted by~$\pi$.
Let~$\pirig\colon \Yrig\rightarrow \Xrig$ be the induced morphism
of rigid~$K$-spaces \'a la Raynaud; cf. \S\ref{subsection:recall}.
In \S\ref{section:measure} we define a ``measure of
singularity"~$\nu_\gerX(P)\in \QQ^{\geq 0}$ (respectively,
~$\nu_\gerY(Q)\in \QQ^{\geq 0}$) of a point~$P$ of~$\Xrig$
(respectively, ~$\Yrig$); the definition is modelled over the
notion of measure of supersingularity for modular curves. For
every interval~$I \subset \RR$ we have an admissible open
set~$\gerY_\rig I$, whose closed point are~$\{Q\in \Yrig:
\nu_\gerY(Q) \in I\}$. The set~$\gerX_\rig I$ is defined
similarly. The following theorem is proven in \S\ref{section:main
theorem}.

\medskip

\id {\bf Theorem A.}\; \emph{Assume~$e>1$. The morphism~$\pirig\colon \Yrig\rightarrow \Xrig$ admits a section
\[\gers_\rig\colon \gerX_\rig[0,e/(e+1)) \rightarrow \Yrig.\]
This section is maximal, namely, it can not be extended to any
connected admissible open properly containing
$\gerX_\rig[0,e/(e+1))$. }
\medskip

\id The reader acquainted with the theory of canonical subgroups
will recognize that this theorem implies the classical existence
theorem for canonical subgroups over modular curves, including the
further statement (that to the best of our knowledge is not
recorded in the literature) that the region over which one defines
the canonical subgroup is the maximal possible, even from the
point of view of maps of rigid spaces. The following theorem,
proven in \S~\ref{section:involution}, will also be familiar to
that reader as giving the behavior of the measure of
supersingularity upon passing to a quotient by a subgroup of
order~$p$. We introduce the following terminology: Let~$Q \in
\Yrig$. We say that~$Q$ is (i) {\it canonical} if~$\nu_\gerY(Q) <
e/(e+1)$; (ii) {\it anti-canonical} if~$\nu_\gerY(Q) > e/(e+1)$;
(iii) {\it too singular} if~$\nu_\gerY(Q)=e/(e+1)$.

\medskip

\id {\bf Theorem B.} \emph{Let~$w$ be an automorphism of~$\gerY$ that
permutes the components of~$\gerY$. We denote by~$w$ also the
induced automorphism of~$\Yrig$ and its effect of points by~$Q
\mapsto Q^w$.
\begin{enumerate}
\item~$\nu_\gerX(\pirig Q) = 0 \;\Leftrightarrow\;
\nu_\gerX(\pirig Q^w) = 0$. In this case~$Q$ is canonical if and
only if~$Q^w$ is anti-canonical. \item If~$\nu_\gerX(\pirig Q)
<(e+1)^{-1}$ and~$Q$ canonical, then~$\nu_\gerX(\pirig Q^w) =
e\cdot \nu_\gerX(\pirig Q)$ and~$Q^w$ is anti-canonical. \item
If~$\nu_\gerX(\pirig Q)=(e+1)^{-1}$ and~$Q$ is canonical
then~$Q^w$ is too singular. \item If~$(e+1)^{-1} <\nu_\gerX(\pirig
Q) < e(e+1)^{-1}$  and~$Q$ is canonical, then~$\nu_\gerX(\pirig
Q^w) = 1 - \nu_\gerX(\pirig Q)$ and~$Q^w$ is canonical. \item
If~$\nu_\gerX(\pirig Q) < e(e+1)^{-1}$ and~$Q$ is anti-canonical,
then~$\nu_\gerX(\pirig Q^w) = e^{-1}\nu_\gerX(\pirig Q)$,
and~$Q^w$ is canonical. \item~If~$Q$ is too singular
then~$\nu_\gerX(\pirig Q^w)=(e+1)^{-1}$ and~$Q^w$ is canonical.
\end{enumerate}}

\

\

\id \emph{Acknowledgments:} The authors benefited from an example
of R. Coleman that inspired the proof of
Proposition~\ref{prop:Coleman}. The first name author was
partially supported by an NSERC grant no. 227040. The second-named
author would like to thank CICMA and the department of mathematics
at McGill university for their support and hospitality.

\

\

\section{Background material}
\subsection{Rigid analytic varieties and formal schemes}\label{subsection:recall}
We recall here the connection between rigid analytic varieties and
formal schemes as developed by Raynaud and Berthelot. Our
exposition follows \cite{BoschLutkebohmertI, BoschLutkebohmertII,
Berthelot, deJong}.

\

\id Let~$R$ be a valuation ring of Krull dimension~$1$, complete
and separated with respect to the~$\gerI$-adic topology, where
$\gerI = (\varpi)$ is contained in the maximal ideal of~$R$. Let
$K$ be the field of fractions of~$R$. For free variables~$\xi =
(\xi_1, \dots, \xi_n)$ we let~$R\langle \xi\rangle = \{ \sum_\nu
c_\nu\xi^\nu \in R[\![\xi]\!]: \lim c_\nu = 0\}$ be the strictly
convergent powerseries, i.e. precisely those that converge on the
polydisc~$\{(a_1, \dots, a_n) : \vert a_i \vert \leq 1\; \forall
i\}$.

Recall that for a general commutative ring~$B$ and an ideal~$J$
of~$B$ one defines the \emph{$J$-torsion} of~$B$ as the ideal~$\{
b\in B: J^n b = 0 \; \text{\rm for some } n\in \NN\}$. If~$J =
(g_1, \dots, g_r)$, the~$J$-torsion is the kernel of the canonical
homomorphism~$R \rightarrow \prod_{i = 1}^r R[g_i^{-1}]$. If this
ideal is~$\{ 0 \}$ we say that~$B$ has no~$J$-torsion.

An \emph{admissible~$R$-algebra} is an~$R$-algebra with no
$\gerI$-torsion (equivalently, flat over~$R$) that is isomorphic
to~$R\langle \xi\rangle/\gera$, where~$\xi = (\xi_1, \dots,
\xi_n)$ for some integer~$n$; it implies that~$\gera$ is a
finitely generated ideal. For us, the admissible~$R$-algebras are
the building blocks of two different categories - a category of
rigid spaces and a category of formal schemes.

An affine formal~$R$-scheme~$\gerX$ is called admissible if it
is of the type~$\gerX = \Spf(A)$, where~$A$ is an admissible
$R$-algebra. We may then write~$\gerX = \underset{\lambda
\rightarrow \infty}{\lim} \gerX_\lambda$, where~$\gerX_\lambda :=
\gerX\otimes (R/(\varpi^\lambda))$, $\lambda \in \NN$, can be
identified with the scheme $\Spec(A\otimes R/(\varpi^\lambda))$.
Being admissible is a local property and so one gets a natural
definition of an \emph{admissible formal~$R$-scheme}.

The notion of admissible blow-up is needed to define an
equivalence of categories between a category of formal schemes and
a category of rigid spaces. The definition of admissible formal
blow-up is designed to be local on the base. We review, thus, only
the affine case. Let~$\gerX = \Spf (A)$ be an affine
admissible~$R$-formal scheme, $A = R\langle \xi\rangle/\gera$.
Let~$\scrA$ be an open ideal, i.e., containing~$(\varpi^\lambda)$
for some~$\lambda>0$. The \emph{admissible formal blow-up}
of~$\gerX$ at~$\scrA$ is~$\gerX^\prime = \underset{\lambda \rightarrow
\infty}{\lim}\; {\rm Proj} \oplus_{n=0}^\infty\left( \scrA^n
\otimes \calO_\gerX/(\varpi^\lambda)\right)$ with the canonical
map~$\varphi\colon \gerX^\prime \rightarrow \gerX$. Then~$\gerX^\prime$ is an
admissible formal~$R$-scheme over which~$\scrA\calO_{\gerX^\prime}$ is
invertible.

Suppose that~$\scrA = (f_0, \dots, f_m)$ and let
$\widetilde{\varphi}\colon \widetilde{\gerX^\prime} \rightarrow
\widetilde{\gerX}=\Spec(A)$ be the usual scheme theory blow-up
of~$A$ at the ideal~$\scrA$. Then, upon
taking~$(\varpi)$-completion of~$\widetilde{\varphi}\colon
\widetilde{\gerX^\prime} \rightarrow \widetilde{\gerX}$ we
get~$\varphi\colon \gerX^\prime \rightarrow \gerX$. On the other hand,
$\widetilde{\varphi}\colon \widetilde{\gerX^\prime} \rightarrow
\widetilde{\gerX}$ admits a local description. The
scheme~$\widetilde{\gerX^\prime}$ has an affine cover by~$\{
\Spec(A_i^\prime): i = 0, 1, \dots, m\}$, where~$A_i^\prime = A_i^{\prime\prime}/(f_i-{\rm
torsion})$ and~$A_i^{\prime\prime} = A\left[\frac{f_0}{f_i}, \dots,
\frac{f_m}{f_i} \right] = A\left[\frac{\zeta_0}{\zeta_i}, \dots,
\frac{\zeta_m}{\zeta_i} \right]/\left( f_i\frac{\zeta_j}{\zeta_i}
- f_j\right)$. Then the~$(\varpi)$-completions of~$A_i^\prime, A_i^{\prime\prime}$
are given by~$\hat{A}_i^\prime = \hat{A}_i^{\prime\prime}/(f_i-{\rm torsion})$ and
$\hat{A}_i^{\prime\prime} = A\left\langle\frac{f_0}{f_i}, \dots,
\frac{f_m}{f_i}\right\rangle = A\left\langle
\frac{\zeta_0}{\zeta_i}, \dots, \frac{\zeta_m}{\zeta_i}
\right\rangle/\left( f_i\frac{\zeta_j}{\zeta_i} - f_j\right)$;
they give rise to an affine covering~$\{ \Spf(\hat{A}_i^\prime):i = 0,
1, \dots, m\}$ of~$\gerX^\prime$.

\

\id For an admissible~$R$-algebra~$A = R\langle \xi\rangle/\gera$,
let~$A_{\rm rig}: = A\otimes_R K = K\langle \xi\rangle/\gera
K\langle \xi\rangle$ -- this is an affinoid~$K$-algebra. This
construction extends to provide a functor \[{\bf rig}\colon
\{\text{admissible formal~$R$-schemes}\} \;\rightarrow\;
\{\text{rigid~$K$-spaces}\},\qquad\gerX \mapsto \gerX_\rig.\] One
calls~$\gerX_\rig~$ the \emph{generic fibre} of the
formal~$R$-scheme~$\gerX$.
\begin{thm}\emph{(Raynaud)} The functor~${\bf rig}$ is an equivalence
of categories between (i) the category of quasi-compact admissible
formal~$R$-schemes, localized by admissible formal blow-ups, and
(ii) the category of quasi-compact and quasi-separated rigid
$K$-spaces.
\end{thm}
\id It is easy to see from the construction that a flat morphism
of formal schemes induces a flat morphism of rigid spaces. The
converse is also true \cite[Thm. 5.2]{BoschLutkebohmertII}: every
flat morphism of rigid~$K$-spaces comes from a flat morphism of
suitable formal schemes yielding the given rigid spaces. A flat
morphism in the category of rigid spaces has image which is a
finite union of affinoids, in particular, it is open
\cite[Cor.~5.11]{BoschLutkebohmertII}.

We will need to use the \emph{specialization map}. In the affine
case, the points of~$\gerX_\rig$ are the maximal ideals of the
algebra~$A\otimes_R K$; these are in bijection with quotients
of~$A$ that are integral, finite and flat over~$R$. If~$T$ is such
a quotient (it is the valuation ring of a finite extension of
$K$), corresponding to a point~$t\in \gerX_\rig$, we get a closed
immersion of formal $R$-schemes~$\Spf(T) \rightarrow \Spf(A)$,
whose image is supported on a closed point of~$\gerX$ that we
denote by~$\spe(t)$. The definition can be extended to any
formal~$R$-scheme. We get a morphism of ringed spaces~$\spe\colon
\gerX_\rig \rightarrow \gerX$ \cite[IV 4.9]{SGA4}. For every
affine open~$U = \Spf(B) \subset \gerX$, we have ~$\spe^{-1}(U) =
U_\rig$.

\

\id Assume that~$R$ is a discrete valuation ring with residue
field~$\kappa$. In \cite{Berthelot} Berthelot generalizes the
above construction to associate a generic fibre to any locally
noetherian formal scheme~$\gerX$ flat over~$R$ that satisfies a
condition weaker than admissibility: that the special fibre
of~$\gerX$, denoted by~$\gerX_0$ and defined by the ideal of
definition~$\scrI$, is a scheme locally of finite type
over~$\kappa$. This condition is independent of the choice
of~$\scrI$ and coincides with admissibility if~$\varpi\calO_\gerX$
is an ideal of definition for~$\gerX$. We will describe the
construction in the affine case. Let~$\gerX=\Spf(A)$
and~$I=H^0(\gerX,\scrI)$ with generators~$g_1,\dots,g_r$. For~$n
\geq 1$ define \[A_n=A\langle T_1,\dots, T_r\rangle/(g_1^n-\varpi
T _1,\dots,g_r^n-\varpi T_r).\] The condition on~$\gerX$ implies
that~$A_n/\varpi A_n$ is finitely generated over~$\kappa$, and
hence~$\gerX^n=\Spf(A_n)$ is an admissible formal scheme over~$R$.
Applying Raynaud's construction we obtain a rigid analytic
space~$\gerX^n_{\rig}$. For~$m>n$ we have a homomorphism~$A_m
\rightarrow A_n$, defined by sending~$T_i$ to~$g_i^{m-n}T_i$,
inducing a morphism of rigid spaces~$\gerX^n_{\rig} \rightarrow
\gerX^m_{\rig}$. It is easy to see that this morphism is an open
immersion and identifies~$\gerX^n_{\rig}$ with the subdomain of
$\gerX^m_{\rig}$ over which~$|g_i(x)|\leq |\varpi|^{1/n}$. The
generic fibre of~$\gerX$, denoted as before by~$\gerX_{\rig}$, is
defined to be the union of~$\gerX^n_{\rig}$ via the above
inclusions. The rigid spaces~$\gerX^n_\rig$ form an admissible
cover of~$\gerX_\rig$. The construction yields a functor~${\bf
rig}$ whose target is the category of quasi-separated rigid
$K$-spaces.

As an illustration, take~$\gerX$ to be
$\Spf(R[\![\xi_1,...,\xi_r]\!])$ with the ideal of definition
$I=(\varpi,\xi_1,...,\xi_r)$. Then~$\gerX_{\rig}$ is simply the
open unit polydisc of dimension~$r$, which is not quasi-compact,
and~$\gerX^n_{\rig} \subset \gerX_{\rig}$ is the affinoid
subdomain over which~$|\xi_i| \leq |\varpi|^{1/n}$, which is
isomorphic to a closed unit polydisc, and hence is quasi-compact.
Similarly, for $\gerX=\Spf(R[\![x_1,x_2]\!]/(x_1x_2-a))$, where $a
\in R$, one sees that $\gerX_\rig$ is the open  annulus over $K$
with radii $(|a|,1)$.

As in the admissible case, one can define a specialization map
$\spe\colon \gerX_{\rig} \rightarrow \gerX$ by taking the direct
limit of the maps~$\gerX^n_{\rig} \overset{\spe}{\longrightarrow}
\gerX^n \rightarrow \gerX$. The following is Proposition 0.2.7. of \cite{Berthelot}.

\begin{prop}\label{prop:Berthelot}
Let notation be as above. Let~$\gerZ \subseteq \gerX_0$ be a
closed subscheme. Let~$\gerX^{\wedge \gerZ}$ denote the formal
completion of~$\gerX$ along~$\gerZ$. Then~$\spe^{-1}(\gerZ)$ is an
admissible open subset of~$\gerX_\rig$ and the canonical
morphism~$\gerX^{\wedge \gerZ}_\rig \rightarrow \gerX_\rig$
induces a functorial isomorphism~$\gerX^{\wedge \gerZ}_\rig \cong
\spe^{-1}(\gerZ)$.
\end{prop}

\

\subsection{Algebraic geometric input}\label{subsection:algeom}
As in the Introduction, let~$\calO$ be the ring of integers of a
finite extension~$K$ of~$\QQ_p$, $\varpi$ a uniformizer of~$\calO$
and~$\kappa = \calO/(\varpi)$ the residue field. Let~$X, Y$ be
relative curves over~$\calO$. We assume that~$X, Y$ are regular
schemes, $X \rightarrow \Spec(\calO)$ is smooth and~$\pi\colon Y
\arr X$ is a finite flat morphism of degree~$e+1$. Moreover, we
assume that (i) there exists a section~$s\colon X\otimes \kappa
\arr Y\otimes \kappa$ to~$\pi \otimes \kappa$, that (ii) the
special fibre~$Y\otimes \kappa$ is a reduced normal crossing
divisor with two components, and that (iii) the set theoretic
preimage~$(\pi\otimes \kappa)^{-1}(\pi\otimes \kappa)(Q)$ is equal
to~$Q$ for any singular point~$Q\in Y\otimes \kappa$.

\

\id The following lemma must be known to the experts; for lack of
a reference we provide a proof.
\begin{lem}\label{lemma:classification of local rings}
Let~$(A, \germ)$ be a regular two-dimensional complete local ring
containing~$\calO$, such that~$\calO$ is integrally closed
in~$A$,~$\germ \cap \calO=(\varpi)$, and $\kappa \subseteq
A/\germ$ is an algebraic extension.
\begin{enumerate}
\item If~$A\otimes \kappa$ is regular then~$A \cong \calO[\![x]\!]$.
\item If ~$A\otimes \kappa\cong \kappa[\![s, t]\!]/(st)$
then~$A \cong \calO[\![x, y]\!]/(xy - \varpi)$.
\end{enumerate}
\end{lem}
\begin{proof} First note that~$A/\germ \supseteq \kappa$ and so the
local homomorphism~$W(A/\germ) \rightarrow A$ has image
containing~$W(\kappa)$ viewed as a subring of~$\calO$.
Since~$A/\germ$ is an algebraic extension of~$\kappa$,
$W(A/\germ)$ is integral over $W(\kappa)$. Since~$\calO$ is
integrally closed in~$A$ it follows that~$W(A/\germ)$ is contained
in~$\calO$. In particular, $A/\germ=\kappa$.

If~$A\otimes \kappa$ is regular it follows by Cohen's Theorem
that~$A\otimes \kappa \cong \kappa[\![x]\!]$. This gives a  morphism
$\calO[\![x]\!] \rightarrow A$ which is surjective by Nakayama's
lemma; since both rings are domains of the same dimension, we
conclude that~$\calO[\![x]\!] \rightarrow A$ is an isomorphism
(the kernel is a prime ideal of height~$0$).

Assume then that~$A\otimes \kappa\cong \kappa[\![s, t]\!]/(st)$.
Let~$x^\prime, y^\prime \in A$ be elements reducing to~$s, t$,
respectively. The homomorphism~$\calO[\![x, y]\!] \rightarrow A$,
taking~$x, y$ to~$x^\prime, y^\prime$ respectively, is surjective
by Nakayama's lemma. Let~$\gerp$ be the kernel; it is a prime
ideal of height~$1$. In fact~$\gerp$ is a principal ideal,
because~$\calO[\![x, y]\!]$ is a factorial ring and by a theorem
of Krull every prime ideal of height~$1$ is principal. We may
therefore write~$\gerp = (h(x, y))$, where~$h(x, y) = xyv - \varpi
z$ for some~$v, z \in \calO[\![x, y]\!]$. It follows
that~$A\otimes \kappa \cong \kappa[\![x, y]\!]/(xy\bar{v})$,
where~$\bar{v}$ is the reduction of~$v$ modulo~$\varpi$.
Since~$\kappa[\![x,y]\!]/(xy\bar{v})\cong \kappa[\![s,t]\!]/(st)$
by the map taking $x\mapsto s$ and $y\mapsto t$, it follows
that~$\bar{v}$ is a unit. This implies that $v$ itself is a unit
and so~$A\cong \calO[\![x, y]\!]/(xy - \varpi z)$.

We next claim that the ring~$A$ is regular if and only if~$z$ is a
unit. Indeed, if $z$ is a unit then $A\cong \calO[\![x,
yz^{-1}]\!]/(x\cdot yz^{-1}-\varpi)$, which is easily checked to
be regular. Assume now that~$A$ is regular. Then~$(\varpi, x,
y)/I$ is a~$2$-dimensional~$\kappa = A/\germ$ vector space,
where~$I = (\varpi, x, y)^2 + (xy - \varpi z)$ and~$\germ$, the
maximal ideal of~$A$, is the image of~$(\varpi, x, y)$. Hence, for
some~$c_1, c_2, c_3 \in A$, not all in~$\germ$, we have~$c_1
\varpi + c_2 x + c_3 y \in I$. Such a relation gives
modulo~$\varpi$ the relation $c_2x + c_3y \in (x, y)^2$. Since the
cotangent space at the singular point is two dimensional with
basis $\{x,y\}$, it follows that modulo~$\varpi$ we have~$c_2, c_3
\in (x,y)$. Thus, we must have~$c_2, c_3 \in \germ$. Therefore,
$A$ is regular implies that~$\varpi\in I$. Thus, $\varpi
\pmod{\germ^2} \in I/\germ^2 = (\varpi z)/\germ^2$. It follows
that~$z$ is a unit modulo~$\germ^2$ and hence is a unit.
\end{proof}

\begin{lem}\label{lem:localcoordinates}Let~$Q\in Y$ be a singular point and~$ P= \pi(Q)$.
There is a choice of local coordinates at~$Q$ and~$P$
giving~$\calO_Y^{\wedge Q}\cong \calO[\![x, y]\!]/(xy - \varpi)$
and~$\calO_X^{\wedge P} \cong \calO[\![t]\!]$, respectively, such
that on the level of completed local rings at~$Q$ and~$P$ the
morphism~$\pi$ is given by
\begin{equation}t\mapsto x + (yu)^e +
f(y) + \varpi g,\end{equation} where~$f(y) \equiv 0
\pmod{(y^{e+1})}$ and~$u$ is a unit congruent to~$1$
modulo~$\varpi$.
\end{lem}
\begin{proof} It follows from Lemma~\ref{lemma:classification of local rings}
that the map~$Y\rightarrow X$ can be written at a singular
point~$Q\in Y$ in the form of an~$\calO$-algebra local
homomorphism~$\pi^*\colon \calO[\![t]\!] \rightarrow \calO[\![x,
y]\!]/(xy - \varpi)$. Now, upon reduction modulo~$\varpi$, we get
a homomorphism of~$\kappa$-algebras~$\pi^*\otimes \kappa\colon
\kappa[\![t]\!] \rightarrow \kappa[\![x, y]\!]/(xy)$. By our
assumptions on~$\pi\otimes \kappa$, the
compositions~$\kappa[\![t]\!] \rightarrow \kappa[\![x, y]\!]/(xy)
\underset{y \mapsto 0}{\rightarrow} \kappa[\![x]\!]$
and~$\kappa[\![t]\!] \rightarrow \kappa[\![x, y]\!]/(xy)
\underset{x \mapsto 0}{\rightarrow} \kappa[\![y]\!]$ are given,
w.l.o.g., by~$t \mapsto x$ and~$t\mapsto y^e+ f_1(y)$,
where~$f_1(y) \equiv 0 \mod{(y^{e+1})}$ (the existence of the
section implies that every ramification index is equal to~$e$).
Thus, the map~$\pi^\ast \otimes \kappa$ is determined by the image
of~$t$ which has the form~$x+y^e + f_1(y) + xyf_2(x, y)$.

Our goal now is to change coordinates on~$A := \calO[\![x,
y]\!]/(xy - \varpi)$ so as to simplify this map and still have the
same presentation, namely, find~$\hat x, \hat y\in A$ such
that~$\calO[\![x, y]\!]/(xy - \varpi) = \calO[\![\hat x, \hat
y]\!]/(\hat x\hat y - \varpi)$. First note that since~$A$
is~$\varpi$-adically complete the map of units~$A^\times
\rightarrow (A\otimes \kappa)^\times$ is surjective. Let~$u^\prime = (1 +
yf_2(x, y)) \in (A\otimes \kappa)^\times$ and~$\hat u$ any lift of
it to~$A^\times$. Let~$\hat x = x\hat u, \hat y = y \hat u^{-1}$.
Then we have~$\calO[\![\hat x, \hat y]\!](\hat x \hat y - \varpi)
= \calO[\![x, y]\!]/(xy-\varpi)$ and the map~$\calO[\![t]\!]
\rightarrow \calO[\![\hat x, \hat y]\!]/(\hat x \hat y - \varpi)$
has the form~$t \mapsto \hat x + (\hat y\hat u)^e + \hat f(\hat y)
+ \varpi \hat g$, where~$\hat f$ is a lift of~$f_1$ satisfying
$\hat f(y) \equiv 0 \pmod{(y^{e+1})}$.
\end{proof}

\

\subsection{A measure of singularity}\label{section:measure}

Let~$\pi\colon Y \rightarrow X$ be a morphism of curves as in
\S\ref{subsection:algeom}. We denote by $\gerX$,~$\gerY$ the
formal schemes obtained from~$X$, $Y$ by completion along their
special fibres. Let~$\beta_1,\dots,\beta_h$ be the singular points
of~$Y$. Let~$\alpha_i=\pi(\beta_i)$ for~$i=1,\dots,h$. Recall that
by assumption the $\alpha_i$'s and $\beta_i$'s are defined over
$\kappa$. Let~$D_{\alpha_i}$ (respectively~$D_{\beta_i}$) denote
the inverse image of~$\alpha_i$ (respectively~$\beta_i$) under the
specialization map~$\spe\colon \Xrig \rightarrow \gerX$
(respectively~$\spe\colon \Yrig \rightarrow \gerY$).

By Proposition~\ref{prop:Berthelot}~$D_{\alpha_i}$ is the rigid
space associated to~$\Spf(\calO_X^{\wedge \alpha_i}) \cong
\Spf(\calO[\![t]\!])$, using Lemma~\ref{lemma:classification of
local rings}. Therefore~$D_{\alpha_i}$ is an open disc of
radius~$1$ with parameter~$t$. This parameter is unique up to~$t
\mapsto t^\prime=t\epsilon+\varpi z$, where~$\epsilon \in
\calO^\times$ and~$z \in \calO[\![t]\!]$. For a general closed
point~$P \in D_{\alpha_i}$ the value~$\val(t(P))$ depends on~$t$,
however if~$\val(t(P))<1$, then~$\val(t(P))=\val(t^\prime(P))$ for
any~$t^\prime$ as above. We abuse notation and define
\[\nu_\gerX(P)=\val(t(P)),\] bearing in mind that this is well
defined only if~$\val(t(P)) <1$.

Similarly, $D_{\beta_i}$ is the rigid space associated to
$\Spf(\calO_X^{\wedge \beta_i}) \cong
\Spf(\calO[\![x,y]\!]/(xy-\varpi))$. Therefore, $D_{\beta_i}$ is
an open annulus of radii~$(|\varpi|,1)$ with parameter~$x$. For
any closed point~$Q$ in~$D_{\beta_i}$, define
\[\nu_\gerY(Q)=\val(x(Q)).\]This definition is independent of the
choice of the parameters if chosen as in Lemma
\ref{lem:localcoordinates}. The reason is that any other such
parameter~$x^\prime$ is of the form~$x^\prime=x\epsilon+\varpi z$,
where~$\epsilon \in \calO^\times$ and~$z \in
\calO[\![x,y]\!]/(xy-\varpi)$, and~$\val(x(Q))<1$.

Let~$\scrZ$ denote the complement in~$\Xrig$ of
$\spe^{-1}(\{\alpha_1,\dots,\alpha_h\})$. For a closed point~$P$
in~$\scrZ$ we define~$\nu_\gerX(P)=0$. By Proposition
\ref{prop:section over the ordinary locus} below, the complement
in~$\Yrig$ of~$\spe^{-1}(\{\beta_1,\dots,\beta_h\})$ has two
connected components, $\scrZ^0=\spe^{-1}\left( (Y\otimes
\kappa)^0-\{\beta_1,\dots,\beta_h\} \right)$, and
$\scrZ^\infty=\spe^{-1}\left( (Y\otimes
\kappa)^\infty-\{\beta_1,\dots,\beta_h\} \right)$. For points
in~$\scrZ^\infty$ we define~$\nu_\gerY$ to be~$0$, and
on~$\scrZ^0$ we define~$\nu_\gerY$ to be~$1$. We refer
to~$\nu_\gerX$ and~$\nu_\gerY$ as {\it measures of singularity}.
For an interval~$I$ of real numbers, we define~$\Xrig I$ to be the
set of points of~$\Xrig$ where~$\nu_\gerX$ belongs to~$I$. For~$U$
an admissible open subset of~$\Xrig$ we set~$UI=U \cap \Xrig I$.
We use a similar notation for~$\Yrig$. We call $\scrZ$ the {\it
ordinary locus} of $\Xrig$ and its complement
$\spe^{-1}(\{\alpha_1,\dots,\alpha_h\})=\Xrig(0,\infty)$ the {\it
singular locus} of $\Xrig$. We have
$\pirig^{-1}(\Xrig(0,\infty))=\Yrig(0,1)=\spe^{-1}(\{\beta_1,\dots,\beta_h\})$
which we call the singular locus of $\Yrig$.

\

\

\section{Main Theorem}\label{section:main theorem}

\id In this section we prove Theorem {\bf A} of the Introduction,
using the same notation. Our strategy is to construct sections
separately on the ordinary locus and the singular locus and glue
them by means of a general principle of rigid geometry. We start
by constructing a section to~$\pirig$ over the ordinary locus
of~$\Xrig$.

\begin{prop}\label{prop:section over the ordinary locus}
The map~$\pirig$ induces an isomorphism between~$\scrZ^\infty$ and
$\scrZ$. Therefore there is a unique
section~$\gers_\rig^\infty\colon \scrZ \rightarrow \Yrig$
to~$\pirig$ whose image is~$\scrZ^\infty$. Furthermore,
both~$\scrZ^\infty$ and~$\scrZ^0$ are connected. If $e>1$, then
any section to $\pirig$ on $\scrZ$ coincides with
$\gers_\rig^\infty$.
\end{prop}

\begin{proof} We show the existence of the section on the level of
the formal schemes. The curves~$(Y \otimes \kappa)^\infty$, $(Y
\otimes \kappa)^0$ are connected reduced affine curves. Let~$U$ be
the open subset of~$\gerY$ equal to the underlying set of~$(Y
\otimes \kappa)^\infty \cup (Y \otimes \kappa)^0$. Then~$U$ is
affine in the formal schemes sense, namely, we have an open
immersion~$\Spf(B) \rightarrow \gerY$ whose set theoretic image is
$U$. Under the specialization map~$\spe\colon \gerY_\rig
\rightarrow \gerY$ we have~$\spe^{-1}(U) = \scrZ^0 \cup
\scrZ^\infty$ and, moreover, $\scrZ^0 \cup \scrZ^\infty = U_\rig$
(cf. the discussion in \S~\ref{subsection:recall}). We conclude
the following: We have a morphism $\Spf(B) \rightarrow \Spf(A)$,
induced by a homomorphism of~$\varpi$-adically
complete~$\calO$-algebras $A \rightarrow B$, that yields the
morphism~$\scrZ^0 \cup \scrZ^\infty \rightarrow \scrZ$ and reduces
to the morphism~$Y\otimes\kappa\setminus\{ \beta_i\}_{i=1}^h
\longrightarrow X\otimes \kappa \setminus\{ \alpha_i\}_{i=1}^h$.
It transpires that~$B \otimes \kappa = (A\otimes \kappa) \oplus
B_1$. Using Hensel's lemma to lift idempotents, we conclude that
we have~$B = A^+ \oplus B_1^+$, with $A^+ \otimes \kappa =
A\otimes \kappa, B_1^+ \otimes \kappa = B_1$. Using that~$A
\rightarrow A^+$ is a finite flat homomorphism reducing to an
isomorphism after~$\otimes \kappa$, we conclude that~$A = A^+$.
This gives the existence of the section $\gers^\infty\colon\Spf(A)
\arr \Spf(B)$, the analytification of which is the desired
section~$\gers_\rig^\infty\colon \scrZ \rightarrow \scrZ^0\cup
\scrZ^\infty$ with image~$\scrZ^\infty$. In particular,
$\scrZ^\infty$, being isomorphic to~$\scrZ$, which is a curve
minus finitely many residue discs, is connected.

 Furthermore, the morphism~$\Spf(B_1^+) \rightarrow \Spf(A)$
is finite flat of degree~$e$. To show~$\scrZ^0$ is connected it is
enough to show that~$\Spf(B^+_1)$ is flat over~$\Spf(\calO)$, and
has a reduced and connected special fibre (see
Remark~\ref{rmk:connected}). But this is clear since
~$\Spec(B_1^+\otimes \kappa) = (Y\otimes \kappa)^0$.

For the final assertion, note that the image of any section to
$\pirig$ on $\scrZ$ must be a connected component of
$\pirig^{-1}(\scrZ)=\scrZ^\infty \cup \scrZ^0$, and hence it must
be either $\scrZ^\infty$ or $\scrZ^0$. The latter can't happen
since $\pirig\colon \scrZ^0 \rightarrow \scrZ$ is $e$-to-$1$ and
$e>1$.
\end{proof}
\begin{rmk} \label{rmk:connected} Let~$\gerB = \Spf(B)$ be an admissible formal
scheme, with associated rigid space~$\gerB_\rig$. It is possible
that~$\gerB_\rig$ is disconnected, yet the underlying topological
space of~$\gerB$ is connected. An example is provided when~$\calO$
is a ramified extension of~$\ZZ_p$ and we let~$B = \calO\langle x,
y, T\rangle/(xy - p, (x+y)T - \varpi)$. The associated rigid space
is a disjoint union of two annuli. The special fibre is three
lines meeting at a single point. Note though that~$B\otimes \kappa
= k[x, y, T]/(xy, (x+y)T)$ in which~$xT$ is nilpotent.

On the other hand, assume~$\gerB$ is an admissible formal scheme
over $\calO$ such that $\gerB_\rig$ is affinoid (in particular
$\gerB_\rig=\Sp(B \otimes_\calO K)$ where
$B=H^0(\gerB,\calO_{\gerB})$). If~$\gerB\otimes \kappa$ is
\emph{reduced}, then the connectedness of~$\gerB$ implies the same
for~$\gerB_\rig$. Indeed, if not, then there is a non-trivial
idempotent element $e\in B \otimes_\calO K$. We show that $e \in
B$. Note that by flatness of $\gerB$ over $\calO$ we know that
~$B\injects B\otimes K$. If~$e \not\in B$, we can write~$e=
f/\varpi^n$, where~$n>0$ is minimal, and~$f\in B$. Then we
have~$f^2=\varpi^nf$. Reducing modulo $\varpi$, we
get~$\bar{f}\neq 0$ and~$\bar{f}^2=0$ which contradicts our
assumption on~$\gerB \otimes \kappa$. Therefore~$e\in B$. It then
follows that the decomposition of the ``generic fibre"~$B\otimes
K$, namely of the rigid space, induces a decomposition of the
formal scheme~$\Spf(B)$.
\end{rmk}

\id Let~$\calX$ be a~$K$-rigid analytic space, and~$\calU \subset
\calX$ be an admissible affinoid subdomain. An affinoid subdomain
$\calU \subset \calV \subset \calX$ is called a {\it strict
neighborhood} of~$\calU$ in~$\calX$ if the reduction of the
inclusion~$\iota\colon \calU \rightarrow \calV$ factors through an
affine scheme which is finite over~$\Spec(\kappa)$. See
\cite[\S3]{CGJ} for more details. Any strict neighborhood
of~$\scrZ$ in $\Xrig$ contains a domain of the form~$\Xrig[0,a]$
for some positive~$a \in \QQ$; cf. \cite[Prop. 2.3.2]{KisinLai}.
The following is Lemma 6 of \cite{CGJ}. See also \cite{Berthelot}.

\begin{lem}
Let~$f\colon {\mathcal Y} \rightarrow {\mathcal X}$ be a finite
flat morphism of rigid analytic curves. Let~${\mathcal U}$ be an
affinoid subdomain of~${\mathcal X}$, and~$s\colon {\mathcal U}
\rightarrow {\mathcal Y}$ a section to~$f$. Then~$s$ can be
extended to a strict neighborhood of~${\mathcal U}$ in~${\mathcal
X}$.
\end{lem}

\begin{cor}\label{cor:overconvergence}
The section~$\gers_\rig^\infty$ extends to a section
$\gers_\rig^\dagger$ over~$\Xrig[0,a]$ for some positive~$a \in
\QQ$.
\end{cor}
Next we discuss sections to~$\pirig$ over the singular locus, i.e.
where~$\nu_\gerX >0$.

\begin{prop} \label{prop:extending the section}The map~$\pirig\colon \Yrig(0,1) \rightarrow
\Xrig(0,\infty)$ admits a section~$\gert$ on~$\Xrig(0,e/(e+1))$
whose image is~$\Yrig(0,e/(e+1))$. Such a section is unique. If
$e>1$, then we have the following stronger uniqueness result: any
section to $\pirig$ on a connected admissible open subset of
$\Xrig(0,e/(e+1))$ which contains some circle $D_{\alpha_i}[a,a]$
is obtained by the restriction of $\gert$.
\end{prop}

\begin{proof} We have~$\Yrig(0,1)= \coprod_i D_{\beta_i}$
and~$\Xrig(0,\infty)=\coprod_i D_{\alpha_i}$. Since by our
assumptions~$\pi^{-1}(\alpha_i)=\{\beta_i\}$ as sets, we have
$\pirig^{-1}(D_{\alpha_i})=D_{\beta_i}$, and hence, for the first assertion, it suffices to
show that for each~$i$ the map~$\pirig\colon D_{\beta_i}
\rightarrow D_{\alpha_i}$ admits a section
on~$D_{\alpha_i}(0,e/(e+1))$ whose image
is~$D_{\beta_i}(0,e/(e+1))$. The map~$\pirig\colon D_{\beta_i}
\rightarrow D_{\alpha_i}$ is the analytification of the map
$\pi\colon \Spf(\calO_Y^{\wedge \beta_i}) \rightarrow
\Spf(\calO_X^{\wedge \alpha_i })$ by Proposition
\ref{prop:Berthelot}. By Lemma~\ref{lem:localcoordinates},
choosing local coordinates, this map is given by
\[\calO[\![t]\!] \rightarrow \calO[\![x, y]\!]/(xy - \varpi), \qquad
t \mapsto x + uy^e + f(y) + \varpi g,\] where~$f(y)
\equiv 0 \pmod{y^e}$, $u,g\in \calO[\![x, y]\!]/(xy -\varpi)$,
and~$u$ is a unit. Let~$\tilde{u}, \tilde{g}$ denote arbitrary
liftings of~$u,g$ to~$\calO[\![x,y]\!]$ and define
$g_0(x)=\tilde{g}(x,\varpi/x)$, $u_0(x)=\tilde{u}(x,\varpi/x)$,
and~$f_0(x)=f(\varpi/x)$. Then the map~$\pirig\colon D_{\beta_i}
\rightarrow D_{\alpha_i}$ is the map characterized by \[t(\pirig
Q)=x(Q)+ u_0(x(Q))(\varpi/x(Q))^e+f_0(x(Q))+\varpi g_0(x(Q)).\]

\newpage
\begin{lem} \label{lem: val and pi rig} Let~$Q \in \Yrig$.
\begin{enumerate}
\item If~$\nu_\gerY(Q)<e/(e+1)$ then~$\nu_\gerX(\pirig
Q)=\nu_\gerY(Q)$. \item If~$\nu_\gerY(Q)>e/(e+1)$ then
$\nu_\gerX(\pirig Q)=e(1-\nu_\gerY(Q))<e/(e+1)$. \item If
$\nu_\gerY(Q)=e/(e+1)$ then~$\nu_\gerX(\pirig Q)\geq e/(e+1)$.
\end{enumerate}
\end{lem}

\begin{figure}[ht]
\psfig{figure= 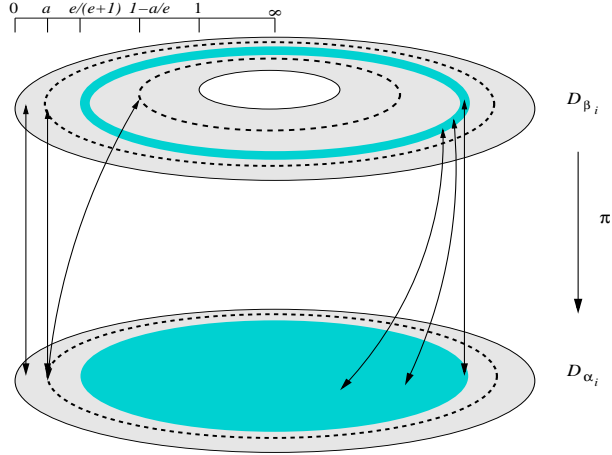,width = 8cm,height = 6cm,
angle=0 }
\caption{The effect of $\pi$ on measures of
singularity.}\label{figureA}
\end{figure}

\begin{proof} The statement is clear for~$Q \in
\scrZ^\infty \cup \scrZ^0$. If~$Q \in D_{\beta_i}$
satisfies~$\nu_\gerY(Q)=\val(x(Q)) < e/(e+1)$, then
\[\val(x(Q))<\min\left\{ \val((\varpi/x(Q))^e ),\;
\val (f_0(x(Q))),\; \val (\varpi g_0(x(Q)))\right\}.\]This implies
that~$\val(t(\pirig Q))=\val(x(Q))$. The other cases are similar.
\end{proof}

From the lemma it follows that \[\pirig^{-1} D_{\alpha_i}(0,e/(e+1))=
D_{\beta_i}(0,e/(e+1)) \coprod D_{\beta_i}(e/(e+1),1).\]
Indeed the lemma proves something stronger: for any~$a \in \QQ$
satisfying~$0 < a < e/(e+1)$ we have

\begin{equation}\label{eqn:radius}
\pirig^{-1}(D_{\alpha_i}[a,a])= D_{\beta_i}[1-a/e,1-a/e] \cup
D_{\beta_i} [a,a].
\end{equation}
This shows that the inverse image of~$D_{\alpha_i}(0,e/(e+1))$
under~$\pirig$ has two connected components each of which maps
onto~$D_{\alpha_i}(0,e/(e+1))$ in a finite flat manner.

We next show that the finite flat morphism~$\pirig\colon
D_{\beta_i}(0,e/(e+1)) \rightarrow D_{\alpha_i}(0,e/(e+1))$ is of
degree one and hence is an isomorphism. The inverse of this map
provides the desired section~$\gert$. To calculate the degree
we restrict the map to a circle~$D_{\alpha_i}[a,a]$ with
$0<a<e/(e+1)$. It is therefore enough to show that~$\pirig\colon
D_{\beta_i}[a,a] \rightarrow D_{\alpha_i}[a,a]$ has degree one. We
show this by reduction modulo~$\varpi$. Our argument is based on
the following general principle.

Let~$\phi\colon \Sp(B) \rightarrow \Sp(A)$ be a finite flat
morphism of~$K$-affinoids. Let~$L$ be a finite field extension
of~$K$ and let~$\phi_L\colon \Sp(B\otimes_K L) \rightarrow
\Sp(A\otimes_K L)$ be the induced morphism. Let~$\theta$ be a
uniformizer of~$L$ and let~$n$ be a positive integer;
let~$(B\otimes_K L)^\circ$ denote the~$\calO_L$-algebra of
functions of supremum norm at most~$1$.
Define~$\overline{B}=(B\otimes_K L)^\circ/(\theta^n)$, and
similarly for~$A$. Let~$\overline{\phi}_L\colon
\Spec(\overline{B}) \rightarrow \Spec(\overline{A})$ be the
induced map. Then, if~$\overline{\phi}_L$ is an isomorphism so
is~$\phi$. The argument reduces to proving that~$\phi_L^*\colon (A
\otimes_K L)^\circ \rightarrow (B \otimes_K L)^\circ$ is
surjective, which, in turn, follows from Nakayama's lemma.

To prove that the reduction of~$\pirig\colon D_{\beta_i}[a,a]
\rightarrow D_{\alpha_i}[a,a]$ is an isomorphism, we first
re-scale. We pass to a finite extension~$L$ of~$K$ with
uniformizer $\theta$ in which there exists an element~$\lambda$ of
valuation~$a$. Setting $x=\lambda x_0$ and~$t=\lambda t_0$ the
map~$\pirig$ becomes a map between circles of radius one
characterized by
\[t_0(\pirig Q)= x_0(Q) + u_0(\lambda x_0(Q))(\varpi^e/\lambda^{e+1})
x_0(Q)^{-e}+\lambda^{-1}f_0(\lambda x_0(Q))+\lambda^{-1}\varpi
g_0(\lambda x_0(Q)).\] Using~$0<a<e(e+1)^{-1}$ and~$f(y)\equiv 0 \pmod{y^{e+1}}$, one
sees that this map reduces modulo~$\theta$ to the identity map of
$ \calO_L/(\theta) [T,1/T]$.

For the second statement we argue as follows. Let~$U \supseteq
D_{\alpha_i}[a,a]$ be a connected admissible open
of~$\Xrig(0,e/(e+1))$ over which there is a section $\gert^\prime$
to $\pirig$. Then $U$, being connected, lies entirely within
$D_{\alpha_i}(0,e/(e+1))$. By Lemma \ref{lem: val and pi rig} the
image of $U$ under $\gert^\prime$ is either a subset of
$D_{\beta_i}(0,e/(e+1))$, or a subset of $D_{\beta_i}(e/(e+1),1)$.
In the former case, by the construction of $\gert$, it is clear
that $\gert^\prime=\gert|_U$. In the latter case,
$\gert^\prime(D_{\alpha_i}[a,a])$ is a connected component of
$D_{\beta_i}[1-a/e,1-a/e]$ by Equation~(\ref{eqn:radius}).
However, since $D_{\beta_i}[1-a/e,1-a/e]$ is connected and
$\pirig\colon D_{\beta_i}[1-a/e,1-a/e] \rightarrow
D_{\alpha_i}[a,a]$ is $e$-to-$1$, we find that $e=1$, which
contradicts our assumption.

\end{proof}

\begin{cor}\label{cor:the section}
 Assume~$e>1$. The morphism~$\pirig\colon \Yrig\rightarrow \Xrig$
admits a unique section
\[\gers_\rig\colon \gerX_\rig[0,e/(e+1)) \rightarrow \Yrig\]
which extends~$\gers_{\rig}^\infty$.
\end{cor}

\begin{proof}
By Corollary~\ref{cor:overconvergence},~$\gers_{\rig}^\infty$
extends to a section~$\gers_{\rig}^\dagger$ on~$\Xrig[0,a]$ for
some positive rational number~$a<e/(e+1)$. By the uniqueness
assertion in Proposition \ref{prop:extending the section} we know
that the restriction of~$\gers_{\rig}^\dagger$ to~$\Xrig(0,a]$ is
obtained as the restriction of~$\gert$. This implies
that~$\gers_{\rig}^\infty$ and~$\gert$ glue together to form the
desired unique section.
\end{proof}

\begin{prop}\label{prop:Coleman} Assume $e>1$. The section $\gers_\rig$
constructed in Corollary \ref{cor:the section} is maximal in the
following sense: Let $U$ be a connected affinoid inside
$\Xrig(0,\infty)$ such that it intersects both $\Xrig(0,e/(e+1))$
and $\Xrig[e/(e+1),1)$ nontrivially  . Then there is no section to
$\pirig$ on $U$.
\end{prop}

\begin{proof} Assume there is a section $\gers$ to $\pirig$ on
such $U$. Since~$U$ is connected it lies inside some
$D_{\alpha_i}$, and intersects both $D_{\alpha_i}(0,e/(e+1))$ and
$D_{\alpha_i}[e/(e+1),e/(e+1)]$ nontrivially. By \cite[\S9.7.2,
Thm. 2]{BGR} any connected affinoid of $D_{\alpha_i}$ is the
complement of a union of finitely many disjoint open discs in a
closed disc. A simple calculation using the non-archimedean
property of the norm shows that a closed disc which intersects
both $D_{\alpha_i}(0,e/(e+1))$ and $D_{\alpha_i}[e/(e+1),e/(e+1)]$
nontrivially, must contain all of $D_{\alpha_i}[e/(e+1),e/(e+1)]$.
Therefore, $U$ contains the complement of a union of finitely many
disjoint open discs $V_1, \dots, V_r$ (which we may assume to have
radius $e/(e+1)$) in the circle $D_{\alpha_i}[e/(e+1),e/(e+1)]$.

We first re-scale as in the proof of
Proposition~\ref{prop:extending the section}: let~$\lambda \in L$
be such that~$\val(\lambda)=e/(e+1)$. Setting~$x=\lambda x_0$
and~$t=\lambda t_0$ the
 map~$\pirig\colon D_{\beta_i}[e/(e+1),e/(e+1)] \rightarrow
D_{\alpha_i}[e/(e+1),\infty)$ becomes a map between a circle~$C$
of radius one and the closed unit disc~$D$ characterized by
\[t_0(\pirig Q)= x_0(Q) + u_0(\lambda x_0(Q))(\varpi^e/\lambda^{e+1})
x_0(Q)^{-e}+\lambda^{-1}f_0(\lambda x_0(Q))+\lambda^{-1}\varpi
g_0(\lambda x_0(Q)).\]  The section $\gers$ is defined on $W$ --
the complement in~$D$ of finitely many residue discs: the open
unit disc together with $\lambda^{-1} V_1, \dots, \lambda^{-1}
V_r$. The reduction~$\bar{\gers}$ of~$\gers\colon W \rightarrow
C$, then, gives a map between ~$\AA^1_{\calO_L/(\theta)}$ minus a
finite number of points (with parameter~$\bar{t}_0$), and
$\AA^1_{\calO_L/(\theta)}$ (with parameter~$\bar{x}_0$)
characterized by
\[\bar{t}_0(Q)= \bar{x}_0({\bar\gers} (Q)) + \overline{(\varpi^e/\lambda^{e+1})}
\bar{x}_0({\bar\gers}(Q))^{-e}.\]Here~$\bar{x}_0({\bar\gers}(-))$
is a rational function in~$\bar{t}_0$ and
$\overline{\varpi^e/\lambda^{e+1}}$ is nonzero by our choice of
$\lambda$. Degree considerations show that this is impossible.
\end{proof}

We summarize the above results as a theorem (Theorem A \, of the
Introduction).

\begin{thm}
Assume~$e>1$. The morphism~$\pirig\colon \Yrig\rightarrow \Xrig$
admits a section
\[\gers_\rig\colon \gerX_\rig[0,e/(e+1)) \rightarrow \Yrig.\]
This section is maximal, namely, it can not be extended to any
connected admissible open properly set containing
$\gerX_\rig[0,e/(e+1))$.
\end{thm}

The canonical subgroup of an elliptic curve can be thought of as a
certain lifting the kernel of Frobenius from characteristic $p$
\cite[Thm. 3.1]{Katz}. We prove a similar result in our setting.
The section $\gers_\rig^\infty=\gers_\rig|_\scrZ$ was constructed
on the level of formal schemes, and by its construction it reduces
to $s$ mod $\varpi$.

Fix $D_{\alpha_i}$ and let $t$ be a coordinate on it, obtained
from an isomorphism~$\calO_X^{\wedge \alpha_i} \cong
\calO[\![t]\!]$ as in Lemma~\ref{lem:localcoordinates}. Also fix
an isomorphism~$\calO_X^{\wedge \beta_i} \cong
\calO[\![x,y]\!]/(xy-\varpi)$ as in loc. cit.; $x$ is a parameter
on~$D_{\beta_i}$ and $xy=\varpi$.

Let $\gamma_P\colon\Sp(L) \rightarrow D_{\alpha_i}$ correspond to
a closed point $P$ which lies in $D_{\alpha_i}(0,e/(e+1))$. Thus
$L$ is a finite extension of~$K$. Let
$\gamma_{\gers_\rig(P)}=\gers_\rig \circ\gamma_P\colon \Sp(L)
\rightarrow D_{\beta_i}$ correspond to the image of $P$ under
$\gers_\rig$. Let $\tilde{\gamma}_P\colon\Spf(\calO_L) \rightarrow
\gerX$ denote the extension of $\gamma_P$ to the formal model, and
similarly define $\tilde{\gamma}_{\gers_\rig(P)}$. Let
$\overline{\gamma}_P$, $\overline{\gamma}_{\gers_\rig(P)}$ denote,
respectively, the reductions of $\tilde{\gamma}_P$,
$\tilde{\gamma}_{\gers_\rig(P)}$ modulo the element $\varpi/t(P)$
of~$\ol$. Let $s^\prime$ denote the base change of $s\colon
X\otimes\kappa \rightarrow Y\otimes\kappa$ from $\kappa$ to
$\calO_L/(\varpi/t(P))$. For simplicity we denote the
$\kappa$-algebra $\calO_L/(\varpi/t(P))$ by $R$.

\begin{prop}
 For closed points $P\in D_{\alpha_i}(0,e/(e+1))$
with $t(P)=r \in \calO_L$ the section $\gers_\rig$ reduces modulo
$\varpi/r$ to $s^\prime$. More precisely, for any $P \in
D_{\alpha_i}(0,e/(e+1))$ we have
$\overline{\gamma}_{\gers_\rig(P)}=s^\prime \circ
\overline{\gamma}_P$.
\end{prop}
\begin{proof} Let us denote the image of an element $a \in \calO_L$ in $R:=\calO_L/(\varpi/t(P))$
by $\overline{a}$. Since $P \in D_{\alpha_i}$, the map
$\tilde{\gamma}_P\colon\Spf(\calO_L) \rightarrow \gerX$ factors
through $\Spf(\calO_X^{\wedge\alpha_i})$. Similarly
$\tilde{\gamma}_{\gers_\rig(P)}$ factors through $\Spf(\calO_Y^{\wedge\beta_i})$. Therefore, it is enough
to prove the statement after replacing $\gerX$ with
$\Spf(\calO_X^{\wedge\alpha_i})$ and $\gerY$ with $\Spf(\calO_Y^{\wedge\beta_i})$.  Then
\[\overline{\gamma}_P\colon\Spec(R) \rightarrow
\Spec(\calO_X^{\wedge\alpha_i} \otimes R)\cong\Spec(R[\![t]\!])\]
is  given by $t \mapsto \overline{t(P)}$. Similarly, the map
\[\overline{\gamma}_{\gers_\rig(P)}\colon\Spec(R) \rightarrow
\Spec(\calO_Y^{\wedge\beta_i} \otimes
R)\cong\Spec(R[\![x,y]\!]/(xy))\] is given by $x \mapsto
\overline{x(\gers_\rig(P))}, y \mapsto
\overline{y(\gers_\rig(P))}$. From the proof of Lemma
\ref{lem:localcoordinates}, we see that the section
\[s^\prime\colon\Spec(R[\![t]\!])\cong\Spec(\calO_X^{\wedge\alpha_i}
\otimes \kappa \otimes_\kappa R) \rightarrow
\Spec(\calO_Y^{\wedge\beta_i}\otimes \kappa \otimes_\kappa R)\cong
\Spec(R[\![x,y]\!]/(xy))\] is given by $x\mapsto t, y\mapsto 0$.
Hence, it is enough to show that $\overline{y(\gers_{\rig}(P))}=0$
and $\overline{x(\gers_{\rig}(P))}=\overline{t(P)}$. For the first
equality notice that by Lemma \ref{lem: val and pi rig} we have
$\val(t(P))=\val(x(\gers_\rig(P)))$ and hence
$y(\gers_\rig(P))=\varpi/x(\gers_\rig(P))$ is divisible by
$\varpi/t(P)$. Since $\val(t(P))=\val(x(\gers_\rig(P)))$, to prove
the second equality it is enough to show that $t(P)$ and
$x(\gers_\rig(P))$ have the same reduction modulo
$y(\gers_\rig(P))=\varpi/x(\gers_\rig(P))$. But that is clear
since from the proof of Proposition \ref{prop:extending the
section} we have
\[t(P)=x(\gers_\rig(P))+ u_0(x(\gers_\rig(P)))(y(\gers_\rig(P)))^e+f_0(x(\gers_\rig(P)))
+\varpi g_0(x(\gers_\rig(P))),\] and  $f_0(x(\gers_\rig(P)))\equiv
0 \pmod {y(\gers_\rig(P))^{e+1}}$.
\end{proof}

\begin{dfn}
Let~$Q \in \Yrig$.

\begin{enumerate}
\item We say that~$Q$ is {\it canonical} if~$Q$ is in the image of
$\gers_\rig$. By the construction of~$\gers_\rig$, this is
equivalent to having~$\nu_\gerY(Q) < e/(e+1)$. If~$Q$ is
canonical, then by Lemma~\ref{lem: val and pi rig} we have
$\nu_\gerY(Q)=\nu_\gerX(\pirig Q)$.

\item We say that~$Q$ is {\it anti-canonical} if~$\nu_\gerY(Q) >
e/(e+1)$. In this case by Lemma~\ref{lem: val and pi rig} we have
$\nu_\gerY(Q)=1-e^{-1}\nu_\gerX(\pirig Q)$.

\item We say that~$Q$ is {\it too singular} if
$\nu_\gerY(Q)=e/(e+1)$. This is equivalent to~$\nu_\gerX(\pirig Q)
\geq e/(e+1)$.
\end{enumerate}
\end{dfn}

\

\

\section{Throwing in an involution}\label{section:involution}

\id In this section we prove the following theorem (Theorem B \,
of the Introduction).

\begin{thm}
Let~$w$ be an automorphism of~$\gerY$ that permutes the components
of~$\gerY$. We denote by~$w$ also the induced automorphism
of~$\Yrig$ and its effect of points by~$Q \mapsto Q^w$. Then:
\begin{enumerate}
\item~$\nu_\gerX(\pirig Q) = 0 \;\Leftrightarrow\;
\nu_\gerX(\pirig Q^w) = 0$. In this case~$Q$ is canonical if and
only if~$Q^w$ is anti-canonical. \item If~$\nu_\gerX(\pirig Q)
<(e+1)^{-1}$ and~$Q$ canonical, then~$\nu_\gerX(\pirig Q^w) =
e\cdot \nu_\gerX(\pirig Q)$ and~$Q^w$ is anti-canonical. \item
If~$\nu_\gerX(\pirig Q)=(e+1)^{-1}$, and~$Q$ is canonical,
then~$Q^w$ is too singular. \item If~$(e+1)^{-1} <
\nu_\gerX(\pirig Q) < e(e+1)^{-1}$, and~$Q$ is canonical,
then~$\nu_\gerX(\pirig Q^w) = 1 - \nu_\gerX(\pirig Q)$ and~$Q^w$
is canonical. \item If~$\nu_\gerX(\pirig Q) < e(e+1)^{-1}$,
and~$Q$ is anti-canonical, then~$\nu_\gerX(\pirig Q^w) =
e^{-1}\nu_\gerX(\pirig Q)$, and~$Q^w$ is canonical. \item~If~$Q$
is too singular, then~$\nu_\gerX(\pirig Q^w)=(e+1)^{-1}$ and~$Q^w$
is canonical.
\end{enumerate}
\end{thm}

We begin by proving the following lemma.

\begin{lem}\label{lem:valuations identity} For any~$Q \in \Yrig$ we have \[\nu_\gerY(Q) +
\nu_\gerY(Q^w)=1.\]
\end{lem}

\begin{proof}
We first note that~$w(\scrZ^\infty)=\scrZ^0$, and hence for~$Q \in
\scrZ^\infty \cup \scrZ^0$ the result follows from the definition
of~$\nu_\gerY$. Assume~$Q \in D_{\beta_i}$ for some~$1 \leq i \leq
h$. The involution~$w$ induces an isomorphism between
$D_{\beta_i}$ and~$D_{\beta_j}$ where~$\beta_j=\beta_i^w$. Let
$x,y$ be coordinates on~$D_{\beta_i}$ as in Lemma
\ref{lem:localcoordinates}. Then~$\eta:=w^*x$ and~$\xi:=w^*y$ are
coordinates on~$D_{\beta_j}$ such that~$D_{\beta_j}$ is the
analytification of~$\Spf(\calO[\![\xi,\eta]\!]/(\xi\eta-\varpi))$.
Because~$w$ switches the two components of~$Y \otimes \kappa$,
$\xi$ is a local parameter on the component containing~$(Y\otimes
\kappa)^\infty$ at the point~$\beta_j$. Examination of the proof
of Lemma~\ref{lem:localcoordinates} shows that there is a local
parameter~$\tau$ on~$D_{\alpha_j}$, and local parameters
$(\hat{\xi}, \hat{\eta})$ on~$D_{\beta_j}$ such that
$\hat{\xi}=\xi\hat{u}$, $\hat{\eta}=\eta\hat{u}^{-1}$, where
$\hat{u}$ is a unit in~$\calO[\![\xi,\eta]\!]/(\xi\eta-\varpi)$,
and such that~$(\tau,\hat{\xi},\hat{\eta})$ are related as in the
statement of Lemma~\ref{lem:localcoordinates}.

By our definition, we can use~$\hat\xi$ to calculate~$\nu_\gerY$
on~$D_{\beta_j}$. Therefore
\[\nu_\gerY(Q^w)=\val(\hat\xi(Q^w))=\val(\xi(Q^w))=\val(y(Q))=1-\val(x(Q))=1-\nu_\gerY(Q).\]
\end{proof}

\id We now prove the theorem.

\medskip

\id (1) is clear.

\medskip

\id (2) Since~$Q$ is canonical, Lemma~\ref{lem: val and pi rig}
implies that~$\nu_\gerY(Q)=\nu_\gerX(\pirig Q)<(e+1)^{-1}$.
Therefore by Lemma~\ref{lem:valuations identity} we have
$\nu_\gerY(Q^w) > e(e+1)^{-1}$, which means that~$Q^w$ is
anti-canonical. It now follows from Lemma~\ref{lem: val and pi
rig} that~$\nu_\gerX(\pirig Q^w) = e(1 - \nu_\gerY(Q^w)) = e
\nu_\gerY(Q) = e \nu_\gerX(\pirig Q)$.

\medskip

\id (3) Since~$Q$ is canonical, Lemma~\ref{lem: val and pi rig}
implies that~$\nu_\gerY(Q)=\nu_\gerX(\pirig Q)=(e+1)^{-1}$, and
therefore~$\nu_\gerY(Q^w)=e(e+1)^{-1}$. This shows that~$Q^w$ is
too singular. It follows from Lemma~\ref{lem: val and pi rig} that
$\nu_\gerX(\pirig Q^w) \geq e(e+1)^{-1}$.

\medskip

\id (4) Since~$Q$ is canonical, we
have~$\nu_\gerY(Q)=\nu_\gerX(\pirig Q)>(e+1)^{-1}$, and hence
$\nu_\gerY(Q^w)<e(e+1)^{-1}$. This shows that~$Q^w$ is canonical.
Therefore, $\nu_\gerX(\pirig Q^w) = \nu_\gerY(Q^w)=1 -
\nu_\gerY(Q)=1-\nu_\gerX(\pirig Q)$.

\medskip

\id (5) Since~$Q$ is anti-canonical, Lemma~\ref{lem: val and pi
rig} shows that ~$\nu_\gerY(Q)=1-e^{-1}\nu_\gerX(\pirig Q) >
e(e+1)^{-1}$. Therefore,
$\nu_\gerY(Q^w)=1-\nu_\gerY(Q)<(e+1)^{-1}$ and hence~$Q^w$ is
canonical. We have~$\nu_\gerX(\pirig Q^w) =
\nu_\gerY(Q^w)=1-\nu_\gerY(Q)=e^{-1}\nu_\gerX(\pirig Q)$.

\medskip

\id(6) Since~$\nu_\gerX(\pirig Q) \geq e(e+1)^{-1}$, by Lemma
\ref{lem: val and pi rig} we have~$\nu_\gerY(Q)=e(e+1)^{-1}$, and
hence~$\nu_\gerY(Q^w)=(e+1)^{-1}$, This shows that~$Q^w$ is
canonical. Therefore, we have~$\nu_\gerX(\pirig
Q^w)=\nu_\gerY(Q^w)=(e+1)^{-1}$.

\

\

\section{Applications}\label{section:applications}

\id In this section we review some of the structure theory for
Shimura curves and show that our results apply to these
situations. Our main references are Drinfeld
\cite{DrinfeldElliptic, Drinfeld} and Carayol \cite{Carayol}. In
particular, we reproduce the classical results on canonical
subgroups \cite{Katz}, as well as more recent developments
\cite{Kassaei1, Kassaei2}.

\

\id Let~$F$ be a totally real field of degree~$d$ with ring of
integers $\calO_F$ and let~$B/F$ be a quaternion algebra split at
exactly one infinite prime of~$F$. Let~$R$ be a maximal order
of~$B$. Let~$\gerp$ be a finite prime of~$F$ at which~$B$
splits,~$F_\gerp$ the completion of~$F$ at the prime~$\gerp$,
$\calO_{F, \gerp}$ its ring of integers with a
uniformizer~$\varpi$, and identify $B\otimes_F F_\gerp$
with~$M_2(F_\gerp)$ so that~$R\otimes_{\calO_F} \calO_{F, \gerp} =
M_2(\calO_{F, \gerp})$. With~$B$ there is associated a projective
system of Shimura curves, initially over the complex numbers but,
by Shimura's theory of canonical models, in fact over~$F$. Let~$G
= {\rm Res}_{F/\QQ}(B^\times)$. Let~$X$ be the~$G(\RR)$-conjugacy
class of the homomorphism~$\CC^\times \rightarrow G(\RR)$
sending~$x+iy$ to $\left[\left(\begin{smallmatrix}x&y \\ - y &
x\end{smallmatrix}\right), 1, \dots, 1\right]\in \GL_2(\RR) \times
(\HH^\times)^{d-1}$. Let~$K$ be an open compact subgroup of
$G(\AA^f)$ of the form $K_\gerp \times K^\gerp$, where~$K_\gerp$
is a subgroup of~$\GL_2(\calO_{F, \gerp})$ and~$K^\gerp$ is ``away
from~$\gerp$". The Shimura curve associated with~$K$ is~$M_K(G,
X)(\CC) = G(\QQ)\backslash G(\AA^f)\times X /K$.

\

\subsection{The case~$F = \QQ$} In this case the Shimura curves
$M_K(G, X)/\QQ$ afford a natural modular description. Consider the
functor associating to a scheme~$S$ the isomorphism classes of
triples~$(A, \iota, \alpha)/S$, where~$A/S$ is an abelian scheme
of relative dimension $2$, $\iota\colon R \rightarrow \End_S(A)$
is an injective ring homomorphism and~$\alpha\colon
\underline{R/NR} \rightarrow A[N]$ is an isomorphism of~$R$-group
schemes; cf. \cite[\S4]{Drinfeld}, \cite[\S4]{DiamondTaylor},
\cite{BuzzardModels}. (Such objects are sometimes called ``false
elliptic curves" because of the similarity with the case of~$B =
M_2(\QQ)$ and the usual modular curves.) This corresponds to the
case where~$K$ is~$\Gamma(N)$ -- the elements of
$(R\otimes_\ZZ\widehat{\ZZ})^\times$ (viewed as a subgroup of
$G(\AA^f)$) that reduce to the identity element under
$(R\otimes_\ZZ\widehat{\ZZ})^\times \rightarrow (R\otimes_\ZZ
\ZZ/N\ZZ)^\times$. For a general~$K$, $K$ contains~$\Gamma(N)$ for
some~$N$ and we take~$\alpha$ up to~$K$-equivalence (\'etale
locally). This makes sense in all characteristics once the level
structure is understood in Drinfeld's sense for which we refer to
\cite{DrinfeldElliptic, KatzMazur}. For~$K$ small enough, there is
therefore a scheme~${\bf M}_K$ over~$\Spec(\ZZ)$ representing this
functor such that~${\bf M}_K \otimes_\ZZ \QQ \cong M_K(G, X)$.

As a module over~$R\otimes \ZZ_p = M_2(\ZZ_p)$, the~$p$-divisible
group~$A[p^\infty]$ of~$A/S$ is a direct sum~$A[p^\infty]_1\oplus
A[p^\infty]_2$ of two isomorphic $p$-divisible groups over $S$,
where the decomposition is determined by the orthogonal
idempotents $\left(\begin{smallmatrix} 1 & 0
\\ 0 & 0 \end{smallmatrix}\right), \left(\begin{smallmatrix} 0 & 0
\\ 0 & 1 \end{smallmatrix}\right) \in M_2(\ZZ_p)$; furthermore, these idempotents are conjugate under
$\left(\begin{smallmatrix} 0 & 1
\\ 1 & 0 \end{smallmatrix}\right)$, which induces the
isomorphism~$A[p^\infty]_1\cong A[p^\infty]_2$. Let~$K^p$ be small
enough and let~$K_p$ be the standard Iwahori subgroup
of~$\GL_2(\ZZ_p)$. The open compact subgroup~$K = K_p \times K^p$
corresponds to a choice of level structure away from~$p$ (given by
$K^p$) and a choice of a non-trivial ideal~$H \subset
M_2(\ZZ/p\ZZ)$. Such~$H$ corresponds, via the $K_p$-equivalence
class of~$\alpha$, to an $R$-invariant subgroup of $A[p]$ of
degree $p^2$. The level structure at $p$ can therefore also be
expressed as an isogeny~$f\colon A_1 \rightarrow A_2$ of false
elliptic curves whose kernel is of degree~$p^2$ and is killed
by~$p$. The conditions on $f$ can also be formulated by requiring
$f$ to have ``false degree" $p$, i.e. that~$f^t \circ f= [p]$ (see
below for the exact meaning of this formula); cf. \cite[p.
453]{DiamondTaylor}, \cite[\S\S10-11]{Kassaei1}.

Let~$A/k$ be a false elliptic curve over an algebraically closed
field~$k$ of characteristic~$p$. One can prove, by means of the
idempotents we have chosen, that the functor of infinitesimal
deformations of~$A$ (resp., together with an Iwahori level
structure $K_p \subset \GL_2(\ZZ_p)$) is equivalent to the functor
of deformation of a~$1$-dimensional~$p$-divisible group of
height~$2$ over~$k$ (resp., with a~$\Gamma_0(p)$-level structure).
Thus, this is exactly the situation arising for elliptic curves
and is well understood; cf. \cite{BuzzardModels}. One concludes
for such choice of $K$ that every geometric connected component of
the special fibre~${\bf M}_K \otimes \FF_p$ of~${\bf M}_K$
consists of two smooth curves crossing transversely at the
supersingular points and so is a normal crossing divisor.
Moreover, the natural morphism~${\bf M}_K \otimes \FF_p \arr {\bf
M}_{\GL_2(\ZZ_p) \times K^p} \otimes \FF_p$ is finite flat of
degree~$p+1$ and admits the usual section taking a false elliptic
curve~$A$ with~$K^p$-structure to~$(A, \Ker(\Fr_A))$ with the
same~$K^p$-structure. The other component is isomorphic to~${\bf
M}_{\GL_2(\ZZ_p) \times K^p} \otimes \FF_p$ as well. Indeed, the
morphism ${\bf M}_K \otimes \FF_p \arr {\bf M}_{\GL_2(\ZZ_p)
\times K^p} \otimes \FF_p$ induces on it a map which is bijective
on geometric points (the pre-image of a point~$A$ is~$(A,
\Ker(\Ver_A)$). Hence the map is purely inseparable of degree~$p$.

There is an involution~$w$ on~${\bf M}_K$ that is best described
by its action on objects: an Iwahori level structure~$f\colon A_1
\rightarrow A_2$ of false elliptic curves is sent by duality
to~$f^t\colon A_2^t \rightarrow A_1^t$. We remark here that every
false elliptic curve carries a principal polarization compatible
with the $R$-action \cite[\S4]{Drinfeld}, hence we get~$f^t\colon
A_2 \rightarrow A_1$, whose isomorphism class is well defined
(independent of the choice of polarization). If the kernel of~$f$
is connected (resp. \'etale) then the kernel of~$f^t$ is \'etale
(resp. connected). It follows that~$w$ permutes the two
irreducible components of every geometric connected component
of~${\bf M}_K \otimes \FF_p$. Finally, there is a finite
extension~$\FF_q \supseteq \FF_p$ over which all the connected
components of ${\bf M}_K\otimes \overline{\FF}_p$ and ${\bf
M}_{\GL_2(\ZZ_p) \times K^p}\otimes \overline{\FF}_p$ are defined
and each connected component is a normal crossing divisor. Using
argument as in Remark~\ref{rmk:connected}, and the fact that ${\bf
M}_K\otimes W(\FF_q)$ is flat over $W(\FF_q)$ and has reduced
special fibre, one find that the connected components of ${\bf
M}_K\otimes \FF_q$ (resp. ${\bf M}_{\GL_2(\ZZ_p) \times
K^p}\otimes \FF_q$) are in bijection with the connected components
of the generic fibre. We conclude that each connected
component~$Y$ of ${\bf M}_K\otimes W(\FF_q)$ and its
image~$X\subseteq {\bf M}_{\GL_2(\ZZ_p) \times K^p}\otimes
W(\FF_q)$ satisfy the hypotheses of this paper. Moreover, a
descent argument, using the uniqueness of the section on each
connected component (see Proposition~\ref{prop:extending the
section}), allows one to get a section over ${\bf M}_{\GL_2(\ZZ_p)
\times K^p}\otimes \QQ_p$ defined over~$\QQ_p$. The application of
our results gives a new proof for the existence and other
properties of canonical subgroups of false elliptic curves,
recovering Theorem 11.1 and Lemma 12.5 of \cite{Kassaei1}.

\

\subsection{The case~$[F:\QQ] = d >1$}
In contrast to the previous case, when~$F\neq \QQ$ there is no
natural modular description of the Shimura curves associated to
$B$. Instead, by making an auxiliary choice of a CM field~$L/F$,
one can associate to the algebra~$B\otimes_F L$ another algebraic
group~$G^\prime$ with the same derived group as that of~$G$. The curves
$M_{K^\prime}(G^\prime, X^\prime)/F$ associated to~$G^\prime$ are PEL Shimura curves.
These auxiliary curves play an important role in Carayol's
construction of an integral model for~$M_K(G, X)/F$ over~$\calO =
\calO_{F, \gerp}$, since they are closely related to the Shimura
curves defined by~$G$ \cite[\S4]{Carayol}. Carayol proves that
such a model~${\bf M}_K$ exists, and that there is a
universal~$p$-divisible~$\calO$-module~$\scrG$ of ($\calO$-)
height~$2$ over the projective limit~${\bf M}_\infty$ of ~${\bf
M}_K$ over $K$. This $p$-divisible group is constructed as a
certain ``piece" of the $p$-divisible group of the universal
abelian variety with additional structure existing over (the
projective limit of) the Shimura curves~${\bf M}_{K^\prime}(G^\prime,
X^\prime)/\calO_F$. Note that the~$p$-divisible group~$\scrG$ does not
carry an~$R\otimes \calO$-structure. In a moral sense, this
structure was already used in reducing the height of
the~$p$-divisible~$\calO$-module to~$2$ (this corresponds to
choosing a particular piece of the~$p$-divisible group of the
universal abelian variety over~${\bf M}_{K^\prime}(G^\prime, X^\prime)/\calO_F$ and
is analogous to the process indicated above for~$F = \QQ$). For
details see \cite{Carayol}, in particular \S\S3.3, 6.3. We discuss
this further.

Assume first that~$K=\GL_2(\calO) \times K^\gerp$. Thus, no level
structure is imposed at~$\gerp$. Carayol constructs
a~$p$-divisible group~$\scrG$ over ${\bf M}_\infty$, which is a
$p$-divisible~$\calO$-module of height~$2$. For any geometric
point $x$ of~${\bf M}_K$, there is a way to define the fibre
$\scrG_x$ by lifting $x$ to a geometric point of ${\bf M}_\infty$.
Over a geometric characteristic~$0$ point~$x$ of~${\bf M}_K$ we
have~$\scrG_x \cong (F_\gerp/\calO)^2$. The prime-to-$\gerp$ level
structure plays a somewhat dormant role. For example, Carayol
proves \cite[\S6.6]{Carayol} a ``Serre-Tate theorem" to the effect
that the formal completion of the henselization of~${\bf M}_K$ at
a geometric point~$x$ of its special fibre pro-represents the
functor of infinitesimal deformations for the~$p$-divisible
$\calO$-module~$\scrG_x$. There are two cases:
\begin{enumerate}\item The \emph{ordinary} case, where~$\scrG_x$ is
isomorphic to~$F_\gerp/\calO \oplus (F_\gerp/\calO)^t$,
where~$(-)^t$ denotes the dual~$p$-divisible group; \item The
\emph{supersingular} case where~$\scrG_x$ is the ``unique"
formal~$\calO$-module of dimension~$1$ and height~$2$ \cite[Prop.
1.7]{DrinfeldElliptic}.
\end{enumerate}
The deformation theory was worked out by Drinfeld. One concludes
that in either case the completed local ring is isomorphic
to~$\hat{\calO}^{\rm nr}[\![t]\!]$ and hence that~${\bf M}_K$ is a
regular surface with a smooth special fibre; cf. \cite[Prop. 4.2,
4.5]{DrinfeldElliptic},\cite[App. \S3]{Carayol}.

Carayol also considers the case of level structure~$K_\gerp(n)
\times K^\gerp$, where~$K_\gerp(n)$ is the subgroup consisting of
matrices in~$\GL_2(\calO)$ congruent to~$1$ modulo~$\gerp^n$.
There is a moduli interpretation of a sort  to the ensuing
morphism $\pi\colon{\bf M}_{K_\gerp(n) \times K^\gerp} \rightarrow
{\bf M}_{K_\gerp(0) \times K^\gerp}$; the group
scheme~$\scrG[\gerp^n]$ descends to ${\bf M}_{K_\gerp(n) \times
K^\gerp}$ and is equipped with a Drinfeld full $\gerp^n$-level
structure, namely, a morphism of~$\calO$-group
schemes~$\alpha\colon \underline{(\gerp^{-n}/\calO)}^2 \rightarrow
\scrG[\gerp^n]$, such that the closed subscheme~$\sum_{P\in
(\gerp^{-n}/\calO)^2}\alpha(P)$ is equal to~$\scrG[\gerp^n]$. The
scheme ${\bf M}_{K_\gerp(n) \times K^\gerp}$ is a torsor over
$\Aut(\underline{(\gerp^{-n}/\calO)}^2) \times {\bf M}_{K_\gerp(0)
\times K^\gerp}$ and the morphism $\pi$ is the natural one (in
particular its fibres are principal homogenous spaces for
$\Aut(\underline{(\gerp^{-n}/\calO)}^2)$). Such level structures
were introduced and studied by Drinfeld in \cite[p.
572]{DrinfeldElliptic}, developed more in \cite[\S7,
Appendix]{Carayol}, and studied extensively in \cite{KatzMazur}.
Again Carayol proves a ``Serre-Tate theorem" as to the nature of
the completed local rings \cite[\S7]{Carayol}. He also proves that
the morphism~$\pi$ extends the natural morphism $M_{K_\gerp(n)
\times K^\gerp}(G, X) \rightarrow  M_{K_\gerp(0) \times
K^\gerp}(G, X)$ induced by the inclusion $K_\gerp(n) \times
K^\gerp\injects K_\gerp(0) \times K^\gerp$.

As Carayol remarks \cite[\S0.4]{Carayol}, the construction and
results extend to any choice of level subgroup at~$\gerp$; in
particular, for $K = K_\gerp\times K^\gerp$, where~$K_\gerp$ is
the Iwahori subgroup. The scheme ${\bf M}_{K_\gerp \times
K^\gerp}$ then carries a finite flat group scheme $\scrH$ (\'etale
locally) with a Drinfeld level structure
$\underline{\gerp^{-1}/\calO}\rightarrow \scrH$ such
that~$\sum_{P\in \gerp^{-1}/\calO}\alpha(P)$ is equal to $\scrH$
as a closed subscheme. The following conclusion follows from
Carayol's work: The completion of the henselization of ${\bf
M}_{K_\gerp \times K^\gerp}$ at a geometric characteristic~$p$
point~$x$ is the ring that pro-represents the functor of
infinitesimal deformations of the divisible $\calO$-module
$\scrG_x$ together with an $\calO$-subgroup scheme of order $q =
\vert \calO/\gerp \vert$ killed by $\gerp$. This moduli problem
can also be phrased in a balanced manner. It can be viewed as
deforming a pair of divisible $\calO$-modules of height $2$, say
$\scrG_x, \scrG_x^\prime$, together with an $\calO$-isogeny
$\scrG_x\rightarrow \scrG_x^\prime$ of degree~$q$ whose kernel is
$\gerp$-torsion.

The situation is again very similar to elliptic curves with
$\Gamma_0(p)$-level structure, and in particular the following
holds. The scheme ${\bf M}_{K_\gerp\times K^\gerp}$ is a regular
two dimensional scheme, flat over $\calO_{F, \gerp}$, the
morphism~$\pi$ is finite flat of degree $q+1$ and the nature
of~$\pi$ at every point is completely understood. In particular,
there are two pre-images to every ordinary point of~${\bf
M}_{\GL_2(\calO)\times K^\gerp}$ and~${\bf M}_{K_\gerp \times
K^\gerp}$ is regular at each; there is a unique pre-image~$y$ to
any geometric supersingular point and the completed local ring
of~$y$ is isomorphic to~$\hat{\calO}^{\rm nr}[\![ s, t ]\!]/(st -
\varpi)$. For completeness we sketch an argument below. We remark
that one can also argue using the results in \cite{Carayol}
obtained for full $\gerp$-level structure. However, Carayol uses
an explicit description of the formal $\calO$-module to obtain his
results. Since we do not anticipate such description to be
available (or indeed useful) in higher-dimensional cases, using
Carayol's result will be contrary to our thesis. We therefore
provide an argument that should extend to the more general
situation we have in mind.

\subsubsection{A sample case} Firstly, we quickly recall the technique of local
models in the particular situation of elliptic curves, which
serves as a good sample case for our problem.

The deformation theory of elliptic curves (or abelian varieties)
can be studied as follows. Given a characteristic~$p$ closed
point~$x$ of a moduli space~${\bf M}$ of elliptic curves with
level prime to $p$ and its universal object~$f\colon \scrE
\rightarrow {\bf M}$, choose an open affine neighborhood~$U\ni x$
and a trivialization of~$\HH^1_{\rm dR}(\scrE/U)\cong \calO_U^2$.
The variation of Hodge structure~$R^0f_\ast \Omega_{\scrE/U}
\rightarrow \HH^1_{\rm dR}(\scrE/U)$ provides a morphism~$U \arr
{\bf Grass}$, where~${\bf Grass}$ is the Grassmann scheme of
locally free, locally direct summands of rank~$1$ of~$\calO_U^2$.
One then shows, using the crystalline theory developed by
Grothendieck, that this morphism is \'etale and so is an
isomorphism on the level of completed local rings of~$x$ and its
image in ${\bf Grass}$; cf. \cite{deJongGamma0(p), DP}. If one
wants to work instead with the $p$-divisible groups, one may
replace~$\HH^1_{\rm dR}$ by a similar object provided by the
theory of displays as developed by Zink, or by the theory of
Cartier-Dieudonn\'e modules, or any other theory studying
deformations of~$p$-divisible groups. For example,
\cite[\S3]{RapoportZink} choose the Lie algebra of the universal
vectorial extension of the $p$-divisible group. By analyzing the
Grassmann scheme, one therefore establishes that the completed
local ring is~$D=W(k(x))[\![t]\!]$.

Under this method, the formal scheme representing the
infinitesimal deformation problem of an elliptic curve with a
subgroup of order~$p$ may be translated to a (formal) incidence
variety. We think of the moduli problem as the one for a cyclic
isogeny~$h\colon  E_1 \rightarrow E_2$ of degree~$p$ between
elliptic curves and we are interested in the completed local ring
of the point on the moduli space that corresponds to such data
over a finite field~$k$ of characteristic~$p$. One may choose the
trivialization of the two $\HH^1_{\rm dR}(\scrE/\Spf(D_i))$, $D_i$
($\cong D$) the completed local ring at $E_i$, such that the
isogeny is given by~$\left(
\begin{smallmatrix} 1 & 0 \\ 0 & p
\end{smallmatrix}\right)$ \cite[\S5.3 ff.]{DP} or \cite{deJongGamma0(p)}.
We are then parameterizing a pair of locally
free,
locally direct summands~$(L_1, L_2)$ of rank~$1$ of~$D^2$ such
that~$\left(
\begin{smallmatrix} 1 & 0 \\ 0 & p
\end{smallmatrix}\right)L_1 \subseteq L_2$. In the ordinary case we
get an $L_1$ whose reduction modulo~$p$ is not killed by $\left(
\begin{smallmatrix} 1 & 0 \\ 0 & 0
\end{smallmatrix}\right)$ and the deformation problem is
represented by the completion of the local ring of a $k$-point~$x$
of~$\PP^1_{W(\FF_p)}$ and so is isomorphic to $W(k)[\![t]\!]$. In
the supersingular case we get an $L_1$ whose reduction
\underline{is} killed by $\left(
\begin{smallmatrix} 1 & 0 \\ 0 & 0
\end{smallmatrix}\right)$. Let $x$ be a $k$-rational point of~$\PP^1_{W(k)}$ and let
$\PP$ be the blow-up of~$\PP^1_{W(k)}$ at~$x$. Its special fibre
has a unique singular point that we shall still denote by~$x$. The
deformation problem is pro-represented by the completion of the
local ring of~$x$ on~$\PP$ and so is isomorphic to
$W(k)[\![s,t]\!]/(st-p)$.

\subsubsection{The calculation of the completed local rings}
Recall that the moduli problem is phrased in a balanced manner.
Let~$x^\prime$ be a closed point of~${\bf M}_{K_\gerp\times
K^\gerp}$ with finite residue field~$k$, and let~$x$ be a $\bar
k$-point supported on~$x^\prime$, where~$\bar k$ is an algebraic
closure of~$k$. The situation we have is of two divisible
$\calO$-modules $\scrG_x, \scrG_x^\prime$ of dimension~$1$ and
height~$2$ over~$\bar k$ and an $\calO$-isogeny $h\colon \scrG_x
\rightarrow \scrG_x^\prime$ of degree~$q = \vert \calO
/(\varpi)\vert$, whose kernel is killed by~$\gerp$.

Let $\scrG$  be $\scrG_x$ or $\scrG_x^\prime$. The Lie algebra of
the universal vectorial extension of $\scrG$, which serves as a
substitute for the first de Rham cohomology, is a free
$\bar{k}$-module of dimension~$2$. As mentioned above, the functor
of infinitesimal deformations of $\scrG$ is pro-representable by
$R^u = \hat \calO^{\rm nr}[\![t]\!]$, which carries a universal
object $\scrG^u$. This can also be proven by the same technique of
local models applied to the relative Lie algebra $\underline{{\rm
Lie}}(\scrG^u)$ of $\scrG^u$ and the Lie algebra $\underline{{\rm
Lie}}(V\scrG^u)\cong (R^u)^2$ of its universal vectorial
extension~$V\scrG^u$, which identifies the completed local ring of
$x$ with the completed local ring of the $\bar k$-point, still
called $x$, on the formal Grassmann scheme $(\PP^1_{R^u})^{\wedge
x}$.

The analogue of Lemma~5.5 of \cite{DP} holds. Namely, one can
choose isomorphisms $\underline{{\rm Lie}}(V\scrG_x^u)\cong
(R^u)^2$ and $\underline{{\rm Lie}}(V\scrG_x'^u)\cong (R^u)^2$
such that $h$ is given by the matrix $A=\left( \begin{smallmatrix}
1 & 0 \\ 0 & \varpi \end{smallmatrix}\right)$ (one should use that
the $p$-divisible groups are in fact polarized, are ``special
$\calO$-modules" in Drinfeld's sense and that $h$ is compatible
with the polarizations). Therefore, the completed local ring of
$\bar x$ is isomorphic to the formal incidence variety in
$(\PP^1_{R^u})^{\wedge x}\times (\PP^1_{R^u})^{\wedge x}$ given by
$A$, i.e., by the closed subscheme over which we have
$A\left(\underline{{\rm Lie}}(\scrG_x^u)\right) \subseteq
\underline{{\rm Lie}}(\scrG_x'^u)$.

In the ordinary case we find that the complete local ring is $\hat
\calO^{\rm nr}[\![t]\!]$, and in the supersingular case we find
that it is $\hat \calO^{\rm nr}[\![s,t]\!]/(st-\varpi)$. Finally,
one may conclude that the completed local ring of $x^\prime$ itself is
$\calO \otimes_{W(\kappa)} W(k)[\![t]\!]$ if $x$ is ordinary and
is $\calO \otimes_{W(\kappa)} W(k_1)[\![s,t]\!]/(st-\varpi)$
if~$x$ is supersingular, where $[k_1: k]\leq 2$.

Given these results, it is straightforward to verify that the
connected components of the generic fibre of a suitable unramified
base-change of ${\bf M}_{K^\gerp \times K_\gerp} \rightarrow {\bf
M}_{\GL_2(\calO) \times K^\gerp}$ satisfy the assumptions of this
paper, including the existence of an involution. In particular,
one has a unique (partial) section on each pair of connected
components of the generic fibres; a descent argument allows one to
conclude that the section can already be defined before
base-change.

\begin{rmk} One may, of course, carry the same analysis for the
Shimura curves $M_{K^\prime}(G^\prime,X^\prime)$. If anything, the
analysis is easier, since it is the one underlying Carayol's
results. Hence, the results of this paper apply to these cases as
well.
\end{rmk}

\begin{rmk} As is clear from our discussion, whenever we are in
a situation of curves $Y \rightarrow X$ such that~$Y$ (or the
fibres) parameterizes group schemes, e.g. in the case of usual
modular curves where $Y$ has a $\Gamma_0(p)$-level structure, or
for pairs ${\bf M}_{K_\gerp \times K^\gerp} \rightarrow {\bf
M}_{K_\gerp(0) \times K^\gerp}$ (or the analogous situations for
the groups $G^\prime$), the construction of a section as in this
paper provides one with a group scheme over the region where the
section is defined. In particular our results reprove Theorems 3.1
and 3.10.7 of \cite{Katz}, and Theorem 9.1 of \cite{Kassaei2} on
canonical subgroups of abelian schemes parameterized by
$M_{K^\prime}(G^\prime,X^\prime)$, and in addition provide an
analogue of Theorem 3.10.7 of \cite{Katz} for such canonical
subgroups.
\end{rmk}

\

\

\end{document}